\documentclass{article}
\long\def\ig#1{\relax}
\ig{Thanks to Roberto Minio for this def'n.  Compare the def'n of
\comment in AMSTeX.}

\newcount \coefa
\newcount \coefb
\newcount \coefc
\newcount\tempcounta
\newcount\tempcountb
\newcount\tempcountc
\newcount\tempcountd
\newcount\xext
\newcount\yext
\newcount\xoff
\newcount\yoff
\newcount\gap%
\newcount\arrowtypea
\newcount\arrowtypeb
\newcount\arrowtypec
\newcount\arrowtyped
\newcount\arrowtypee
\newcount\height
\newcount\width
\newcount\xpos
\newcount\ypos
\newcount\run
\newcount\rise
\newcount\arrowlength
\newcount\halflength
\newcount\arrowtype
\newdimen\tempdimen
\newdimen\xlen
\newdimen\ylen
\newsavebox{\tempboxa}%
\newsavebox{\tempboxb}%
\newsavebox{\tempboxc}%

\makeatletter
\def\settypes(#1,#2,#3){\arrowtypea#1 \arrowtypeb#2 \arrowtypec#3}
\def\settoheight#1#2{\setbox\@tempboxa\hbox{#2}#1\ht\@tempboxa\relax}%
\def\settodepth#1#2{\setbox\@tempboxa\hbox{#2}#1\dp\@tempboxa\relax}%
\def\settokens[#1`#2`#3`#4]{%
     \def\tokena{#1}\def\tokenb{#2}\def\tokenc{#3}\def\tokend{#4}}
\def\setsqparms[#1`#2`#3`#4;#5`#6]{%
\arrowtypea #1
\arrowtypeb #2
\arrowtypec #3
\arrowtyped #4
\width #5
\height #6
}
\def\setpos(#1,#2){\xpos=#1 \ypos#2}

\def\bfig{\begin{picture}(\xext,\yext)(\xoff,\yoff)}
\def\efig{\end{picture}}

\def\putbox(#1,#2)#3{\put(#1,#2){\makebox(0,0){$#3$}}}

\def\settriparms[#1`#2`#3;#4]{\settripairparms[#1`#2`#3`1`1;#4]}%

\def\settripairparms[#1`#2`#3`#4`#5;#6]{%
\arrowtypea #1
\arrowtypeb #2
\arrowtypec #3
\arrowtyped #4
\arrowtypee #5
\width #6
\height #6
}

\def\resetparms{\settripairparms[1`1`1`1`1;500]\width 500}

\resetparms

\def\mvector(#1,#2)#3{
\put(0,0){\vector(#1,#2){#3}}%
\put(0,0){\vector(#1,#2){30}}%
}
\def\evector(#1,#2)#3{{
\arrowlength #3
\put(0,0){\vector(#1,#2){\arrowlength}}%
\advance \arrowlength by-30
\put(0,0){\vector(#1,#2){\arrowlength}}%
}}

\def\horsize#1#2{%
\settowidth{\tempdimen}{$#2$}%
#1=\tempdimen
\divide #1 by\unitlength
}

\def\vertsize#1#2{%
\settoheight{\tempdimen}{$#2$}%
#1=\tempdimen
\settodepth{\tempdimen}{$#2$}%
\advance #1 by\tempdimen
\divide #1 by\unitlength
}

\def\vertadjust[#1`#2`#3]{%
\vertsize{\tempcounta}{#1}%
\vertsize{\tempcountb}{#2}%
\ifnum \tempcounta<\tempcountb \tempcounta=\tempcountb \fi
\divide\tempcounta by2
\vertsize{\tempcountb}{#3}%
\ifnum \tempcountb>0 \advance \tempcountb by20 \fi
\ifnum \tempcounta<\tempcountb \tempcounta=\tempcountb \fi
}

\def\horadjust[#1`#2`#3]{%
\horsize{\tempcounta}{#1}%
\horsize{\tempcountb}{#2}%
\ifnum \tempcounta<\tempcountb \tempcounta=\tempcountb \fi
\divide\tempcounta by20
\horsize{\tempcountb}{#3}%
\ifnum \tempcountb>0 \advance \tempcountb by60 \fi
\ifnum \tempcounta<\tempcountb \tempcounta=\tempcountb \fi
}

\ig{ In this procedure, #1 is the paramater that sticks out all the way,
#2 sticks out the least and #3 is a label sticking out half way.  #4 is
the amount of the offset.}

\def\sladjust[#1`#2`#3]#4{%
\tempcountc=#4
\horsize{\tempcounta}{#1}%
\divide \tempcounta by2
\horsize{\tempcountb}{#2}%
\divide \tempcountb by2
\advance \tempcountb by-\tempcountc
\ifnum \tempcounta<\tempcountb \tempcounta=\tempcountb\fi
\divide \tempcountc by2
\horsize{\tempcountb}{#3}%
\advance \tempcountb by-\tempcountc
\ifnum \tempcountb>0 \advance \tempcountb by80\fi
\ifnum \tempcounta<\tempcountb \tempcounta=\tempcountb\fi
\advance\tempcounta by20
}

\def\putvector(#1,#2)(#3,#4)#5#6{{%
\xpos=#1
\ypos=#2
\run=#3
\rise=#4
\arrowlength=#5
\arrowtype=#6
\ifnum \arrowtype<0
    \ifnum \run=0
        \advance \ypos by-\arrowlength
    \else
        \tempcounta \arrowlength
        \multiply \tempcounta by\rise
        \divide \tempcounta by\run
        \ifnum\run>0
            \advance \xpos by\arrowlength
            \advance \ypos by\tempcounta
        \else
            \advance \xpos by-\arrowlength
            \advance \ypos by-\tempcounta
        \fi
    \fi
    \multiply \arrowtype by-1
    \multiply \rise by-1
    \multiply \run by-1
\fi
\ifnum \arrowtype=1
    \put(\xpos,\ypos){\vector(\run,\rise){\arrowlength}}%
\else\ifnum \arrowtype=2
    \put(\xpos,\ypos){\mvector(\run,\rise)\arrowlength}%
\else\ifnum\arrowtype=3
    \put(\xpos,\ypos){\evector(\run,\rise){\arrowlength}}%
\fi\fi\fi
}}

\def\putsplitvector(#1,#2)#3#4{
\xpos #1
\ypos #2
\arrowtype #4
\halflength #3
\arrowlength #3
\gap 140
\advance \halflength by-\gap
\divide \halflength by2
\ifnum \arrowtype=1
    \put(\xpos,\ypos){\line(0,-1){\halflength}}%
    \advance\ypos by-\halflength
    \advance\ypos by-\gap
    \put(\xpos,\ypos){\vector(0,-1){\halflength}}%
\else\ifnum \arrowtype=2
    \put(\xpos,\ypos){\line(0,-1)\halflength}%
    \put(\xpos,\ypos){\vector(0,-1)3}%
    \advance\ypos by-\halflength
    \advance\ypos by-\gap
    \put(\xpos,\ypos){\vector(0,-1){\halflength}}%
\else\ifnum\arrowtype=3
    \put(\xpos,\ypos){\line(0,-1)\halflength}%
    \advance\ypos by-\halflength
    \advance\ypos by-\gap
    \put(\xpos,\ypos){\evector(0,-1){\halflength}}%
\else\ifnum \arrowtype=-1
    \advance \ypos by-\arrowlength
    \put(\xpos,\ypos){\line(0,1){\halflength}}%
    \advance\ypos by\halflength
    \advance\ypos by\gap
    \put(\xpos,\ypos){\vector(0,1){\halflength}}%
\else\ifnum \arrowtype=-2
    \advance \ypos by-\arrowlength
    \put(\xpos,\ypos){\line(0,1)\halflength}%
    \put(\xpos,\ypos){\vector(0,1)3}%
    \advance\ypos by\halflength
    \advance\ypos by\gap
    \put(\xpos,\ypos){\vector(0,1){\halflength}}%
\else\ifnum\arrowtype=-3
    \advance \ypos by-\arrowlength
    \put(\xpos,\ypos){\line(0,1)\halflength}%
    \advance\ypos by\halflength
    \advance\ypos by\gap
    \put(\xpos,\ypos){\evector(0,1){\halflength}}%
\fi\fi\fi\fi\fi\fi
}

\def\putmorphism(#1)(#2,#3)[#4`#5`#6]#7#8#9{{%
\run #2
\rise #3
\ifnum\rise=0
  \puthmorphism(#1)[#4`#5`#6]{#7}{#8}{#9}%
\else\ifnum\run=0
  \putvmorphism(#1)[#4`#5`#6]{#7}{#8}{#9}%
\else
\setpos(#1)%
\arrowlength #7
\arrowtype #8
\ifnum\run=0
\else\ifnum\rise=0
\else
\ifnum\run>0
    \coefa=1
\else
   \coefa=-1
\fi
\ifnum\arrowtype>0
   \coefb=0
   \coefc=-1
\else
   \coefb=\coefa
   \coefc=1
   \arrowtype=-\arrowtype
\fi
\width=2
\multiply \width by\run
\divide \width by\rise
\ifnum \width<0  \width=-\width\fi
\advance\width by60
\if l#9 \width=-\width\fi
\putbox(\xpos,\ypos){#4}
{\multiply \coefa by\arrowlength
\advance\xpos by\coefa
\multiply \coefa by\rise
\divide \coefa by\run
\advance \ypos by\coefa
\putbox(\xpos,\ypos){#5} }%
{\multiply \coefa by\arrowlength
\divide \coefa by2
\advance \xpos by\coefa
\advance \xpos by\width
\multiply \coefa by\rise
\divide \coefa by\run
\advance \ypos by\coefa
\if l#9%
   \put(\xpos,\ypos){\makebox(0,0)[r]{$#6$}}%
\else\if r#9%
   \put(\xpos,\ypos){\makebox(0,0)[l]{$#6$}}%
\fi\fi }%
{\multiply \rise by-\coefc
\multiply \run by-\coefc
\multiply \coefb by\arrowlength
\advance \xpos by\coefb
\multiply \coefb by\rise
\divide \coefb by\run
\advance \ypos by\coefb
\multiply \coefc by70
\advance \ypos by\coefc
\multiply \coefc by\run
\divide \coefc by\rise
\advance \xpos by\coefc
\multiply \coefa by140
\multiply \coefa by\run
\divide \coefa by\rise
\advance \arrowlength by\coefa
\ifnum \arrowtype=1
   \put(\xpos,\ypos){\vector(\run,\rise){\arrowlength}}%
\else\ifnum\arrowtype=2
   \put(\xpos,\ypos){\mvector(\run,\rise){\arrowlength}}%
\else\ifnum\arrowtype=3
   \put(\xpos,\ypos){\evector(\run,\rise){\arrowlength}}%
\fi\fi\fi}\fi\fi\fi\fi}}

\def\puthmorphism(#1,#2)[#3`#4`#5]#6#7#8{{%
\xpos #1
\ypos #2
\width #6
\arrowlength #6
\putbox(\xpos,\ypos){#3\vphantom{#4}}%
{\advance \xpos by\arrowlength
\putbox(\xpos,\ypos){\vphantom{#3}#4}}%
\horsize{\tempcounta}{#3}%
\horsize{\tempcountb}{#4}%
\divide \tempcounta by2
\divide \tempcountb by2
\advance \tempcounta by30
\advance \tempcountb by30
\advance \xpos by\tempcounta
\advance \arrowlength by-\tempcounta
\advance \arrowlength by-\tempcountb
\putvector(\xpos,\ypos)(1,0){\arrowlength}{#7}%
\divide \arrowlength by2
\advance \xpos by\arrowlength
\vertsize{\tempcounta}{#5}%
\divide\tempcounta by2
\advance \tempcounta by20
\if a#8 %
   \advance \ypos by\tempcounta
   \putbox(\xpos,\ypos){#5}%
\else
   \advance \ypos by-\tempcounta
   \putbox(\xpos,\ypos){#5}%
\fi}}

\def\putvmorphism(#1,#2)[#3`#4`#5]#6#7#8{{%
\xpos #1
\ypos #2
\arrowlength #6
\arrowtype #7
\settowidth{\xlen}{$#5$}%
\putbox(\xpos,\ypos){#3}%
{\advance \ypos by-\arrowlength
\putbox(\xpos,\ypos){#4}}%
{\advance\arrowlength by-140
\advance \ypos by-70
\ifdim\xlen>0pt
   \if m#8%
      \putsplitvector(\xpos,\ypos){\arrowlength}{\arrowtype}%
   \else
      \putvector(\xpos,\ypos)(0,-1){\arrowlength}{\arrowtype}%
   \fi
\else
   \putvector(\xpos,\ypos)(0,-1){\arrowlength}{\arrowtype}%
\fi}%
\ifdim\xlen>0pt
   \divide \arrowlength by2
   \advance\ypos by-\arrowlength
   \if l#8%
      \advance \xpos by-40
      \put(\xpos,\ypos){\makebox(0,0)[r]{$#5$}}%
   \else\if r#8%
      \advance \xpos by40
      \put(\xpos,\ypos){\makebox(0,0)[l]{$#5$}}%
   \else
      \putbox(\xpos,\ypos){#5}%
   \fi\fi
\fi
}}

\def\topadjust[#1`#2`#3]{%
\yoff=10
\vertadjust[#1`#2`{#3}]%
\advance \yext by\tempcounta
\advance \yext by 10
}
\def\botadjust[#1`#2`#3]{%
\vertadjust[#1`#2`{#3}]%
\advance \yext by\tempcounta
\advance \yoff by-\tempcounta
}
\def\leftadjust[#1`#2`#3]{%
\xoff=0
\horadjust[#1`#2`{#3}]%
\advance \xext by\tempcounta
\advance \xoff by-\tempcounta
}
\def\rightadjust[#1`#2`#3]{%
\horadjust[#1`#2`{#3}]%
\advance \xext by\tempcounta
}
\def\rightsladjust[#1`#2`#3]{%
\sladjust[#1`#2`{#3}]{\width}%
\advance \xext by\tempcounta
}
\def\leftsladjust[#1`#2`#3]{%
\xoff=0
\sladjust[#1`#2`{#3}]{\width}%
\advance \xext by\tempcounta
\advance \xoff by-\tempcounta
}
\def\adjust[#1`#2;#3`#4;#5`#6;#7`#8]{%
\topadjust[#1``{#2}]
\leftadjust[#3``{#4}]
\rightadjust[#5``{#6}]
\botadjust[#7``{#8}]}

\def\putsquarep<#1>(#2)[#3;#4`#5`#6`#7]{{%
\setsqparms[#1]%
\setpos(#2)%
\settokens[#3]%
\puthmorphism(\xpos,\ypos)[\tokenc`\tokend`{#7}]{\width}{\arrowtyped}b%
\advance\ypos by \height
\puthmorphism(\xpos,\ypos)[\tokena`\tokenb`{#4}]{\width}{\arrowtypea}a%
\putvmorphism(\xpos,\ypos)[``{#5}]{\height}{\arrowtypeb}l%
\advance\xpos by \width
\putvmorphism(\xpos,\ypos)[``{#6}]{\height}{\arrowtypec}r%
}}

\def\putsquare{\@ifnextchar <{\putsquarep}{\putsquarep%
   <\arrowtypea`\arrowtypeb`\arrowtypec`\arrowtyped;\width`\height>}}
\def\square{\@ifnextchar< {\squarep}{\squarep
   <\arrowtypea`\arrowtypeb`\arrowtypec`\arrowtyped;\width`\height>}}
\def\squarep<#1>[#2`#3`#4`#5;#6`#7`#8`#9]{{
\setsqparms[#1]
\xext=\width                                          
\yext=\height                                         
\topadjust[#2`#3`{#6}]
\botadjust[#4`#5`{#9}]
\leftadjust[#2`#4`{#7}]
\rightadjust[#3`#5`{#8}]
\begin{picture}(\xext,\yext)(\xoff,\yoff)
\putsquarep<\arrowtypea`\arrowtypeb`\arrowtypec`\arrowtyped;\width`\height>%
(0,0)[#2`#3`#4`#5;#6`#7`#8`{#9}]%
\end{picture}%
}}

\def\putptrianglep<#1>(#2,#3)[#4`#5`#6;#7`#8`#9]{{%
\settriparms[#1]%
\xpos=#2 \ypos=#3
\advance\ypos by \height
\puthmorphism(\xpos,\ypos)[#4`#5`{#7}]{\height}{\arrowtypea}a%
\putvmorphism(\xpos,\ypos)[`#6`{#8}]{\height}{\arrowtypeb}l%
\advance\xpos by\height
\putmorphism(\xpos,\ypos)(-1,-1)[``{#9}]{\height}{\arrowtypec}r%
}}

\def\putptriangle{\@ifnextchar <{\putptrianglep}{\putptrianglep
   <\arrowtypea`\arrowtypeb`\arrowtypec;\height>}}
\def\ptriangle{\@ifnextchar <{\ptrianglep}{\ptrianglep
   <\arrowtypea`\arrowtypeb`\arrowtypec;\height>}}

\def\ptrianglep<#1>[#2`#3`#4;#5`#6`#7]{{
\settriparms[#1]%
\width=\height                         
\xext=\width                           
\yext=\width                           
\topadjust[#2`#3`{#5}]
\botadjust[#3``]
\leftadjust[#2`#4`{#6}]
\rightsladjust[#3`#4`{#7}]
\begin{picture}(\xext,\yext)(\xoff,\yoff)
\putptrianglep<\arrowtypea`\arrowtypeb`\arrowtypec;\height>%
(0,0)[#2`#3`#4;#5`#6`{#7}]%
\end{picture}%
}}

\def\putqtrianglep<#1>(#2,#3)[#4`#5`#6;#7`#8`#9]{{%
\settriparms[#1]%
\xpos=#2 \ypos=#3
\advance\ypos by\height
\puthmorphism(\xpos,\ypos)[#4`#5`{#7}]{\height}{\arrowtypea}a%
\putmorphism(\xpos,\ypos)(1,-1)[``{#8}]{\height}{\arrowtypeb}l%
\advance\xpos by\height
\putvmorphism(\xpos,\ypos)[`#6`{#9}]{\height}{\arrowtypec}r%
}}

\def\putqtriangle{\@ifnextchar <{\putqtrianglep}{\putqtrianglep
   <\arrowtypea`\arrowtypeb`\arrowtypec;\height>}}
\def\qtriangle{\@ifnextchar <{\qtrianglep}{\qtrianglep
   <\arrowtypea`\arrowtypeb`\arrowtypec;\height>}}

\def\qtrianglep<#1>[#2`#3`#4;#5`#6`#7]{{
\settriparms[#1]
\width=\height                         
\xext=\width                           
\yext=\height                          
\topadjust[#2`#3`{#5}]
\botadjust[#4``]
\leftsladjust[#2`#4`{#6}]
\rightadjust[#3`#4`{#7}]
\begin{picture}(\xext,\yext)(\xoff,\yoff)
\putqtrianglep<\arrowtypea`\arrowtypeb`\arrowtypec;\height>%
(0,0)[#2`#3`#4;#5`#6`{#7}]%
\end{picture}%
}}

\def\putdtrianglep<#1>(#2,#3)[#4`#5`#6;#7`#8`#9]{{%
\settriparms[#1]%
\xpos=#2 \ypos=#3
\puthmorphism(\xpos,\ypos)[#5`#6`{#9}]{\height}{\arrowtypec}b%
\advance\xpos by \height \advance\ypos by\height
\putmorphism(\xpos,\ypos)(-1,-1)[``{#7}]{\height}{\arrowtypea}l%
\putvmorphism(\xpos,\ypos)[#4``{#8}]{\height}{\arrowtypeb}r%
}}

\def\putdtriangle{\@ifnextchar <{\putdtrianglep}{\putdtrianglep
   <\arrowtypea`\arrowtypeb`\arrowtypec;\height>}}
\def\dtriangle{\@ifnextchar <{\dtrianglep}{\dtrianglep
   <\arrowtypea`\arrowtypeb`\arrowtypec;\height>}}

\def\dtrianglep<#1>[#2`#3`#4;#5`#6`#7]{{
\settriparms[#1]
\width=\height                         
\xext=\width                           
\yext=\height                          
\topadjust[#2``]
\botadjust[#3`#4`{#7}]
\leftsladjust[#3`#2`{#5}]
\rightadjust[#2`#4`{#6}]
\begin{picture}(\xext,\yext)(\xoff,\yoff)
\putdtrianglep<\arrowtypea`\arrowtypeb`\arrowtypec;\height>%
(0,0)[#2`#3`#4;#5`#6`{#7}]%
\end{picture}%
}}

\def\putbtrianglep<#1>(#2,#3)[#4`#5`#6;#7`#8`#9]{{%
\settriparms[#1]%
\xpos=#2 \ypos=#3
\puthmorphism(\xpos,\ypos)[#5`#6`{#9}]{\height}{\arrowtypec}b%
\advance\ypos by\height
\putmorphism(\xpos,\ypos)(1,-1)[``{#8}]{\height}{\arrowtypeb}r%
\putvmorphism(\xpos,\ypos)[#4``{#7}]{\height}{\arrowtypea}l%
}}

\def\putbtriangle{\@ifnextchar <{\putbtrianglep}{\putbtrianglep
   <\arrowtypea`\arrowtypeb`\arrowtypec;\height>}}
\def\btriangle{\@ifnextchar <{\btrianglep}{\btrianglep
   <\arrowtypea`\arrowtypeb`\arrowtypec;\height>}}

\def\btrianglep<#1>[#2`#3`#4;#5`#6`#7]{{
\settriparms[#1]
\width=\height                         
\xext=\width                           
\yext=\height                          
\topadjust[#2``]
\botadjust[#3`#4`{#7}]
\leftadjust[#2`#3`{#5}]
\rightsladjust[#4`#2`{#6}]
\begin{picture}(\xext,\yext)(\xoff,\yoff)
\putbtrianglep<\arrowtypea`\arrowtypeb`\arrowtypec;\height>%
(0,0)[#2`#3`#4;#5`#6`{#7}]%
\end{picture}%
}}

\def\putAtrianglep<#1>(#2,#3)[#4`#5`#6;#7`#8`#9]{{%
\settriparms[#1]%
\xpos=#2 \ypos=#3
{\multiply \height by2
\puthmorphism(\xpos,\ypos)[#5`#6`{#9}]{\height}{\arrowtypec}b}%
\advance\xpos by\height \advance\ypos by\height
\putmorphism(\xpos,\ypos)(-1,-1)[#4``{#7}]{\height}{\arrowtypea}l%
\putmorphism(\xpos,\ypos)(1,-1)[``{#8}]{\height}{\arrowtypeb}r%
}}

\def\putAtriangle{\@ifnextchar <{\putAtrianglep}{\putAtrianglep
   <\arrowtypea`\arrowtypeb`\arrowtypec;\height>}}
\def\Atriangle{\@ifnextchar <{\Atrianglep}{\Atrianglep
   <\arrowtypea`\arrowtypeb`\arrowtypec;\height>}}

\def\Atrianglep<#1>[#2`#3`#4;#5`#6`#7]{{
\settriparms[#1]
\width=\height                         
\xext=\width                           
\yext=\height                          
\topadjust[#2``]
\botadjust[#3`#4`{#7}]
\multiply \xext by2 
\leftsladjust[#3`#2`{#5}]
\rightsladjust[#4`#2`{#6}]
\begin{picture}(\xext,\yext)(\xoff,\yoff)%
\putAtrianglep<\arrowtypea`\arrowtypeb`\arrowtypec;\height>%
(0,0)[#2`#3`#4;#5`#6`{#7}]%
\end{picture}%
}}

\def\putAtrianglepairp<#1>(#2)[#3;#4`#5`#6`#7`#8]{{
\settripairparms[#1]%
\setpos(#2)%
\settokens[#3]%
\puthmorphism(\xpos,\ypos)[\tokenb`\tokenc`{#7}]{\height}{\arrowtyped}b%
\advance\xpos by\height
\advance\ypos by\height
\putmorphism(\xpos,\ypos)(-1,-1)[\tokena``{#4}]{\height}{\arrowtypea}l%
\putvmorphism(\xpos,\ypos)[``{#5}]{\height}{\arrowtypeb}m%
\putmorphism(\xpos,\ypos)(1,-1)[``{#6}]{\height}{\arrowtypec}r%
}}

\def\putAtrianglepair{\@ifnextchar <{\putAtrianglepairp}{\putAtrianglepairp%
   <\arrowtypea`\arrowtypeb`\arrowtypec`\arrowtyped`\arrowtypee;\height>}}
\def\Atrianglepair{\@ifnextchar <{\Atrianglepairp}{\Atrianglepairp%
   <\arrowtypea`\arrowtypeb`\arrowtypec`\arrowtyped`\arrowtypee;\height>}}

\def\Atrianglepairp<#1>[#2;#3`#4`#5`#6`#7]{{%
\settripairparms[#1]%
\settokens[#2]%
\width=\height
\xext=\width
\yext=\height
\topadjust[\tokena``]%
\vertadjust[\tokenb`\tokenc`{#6}]
\tempcountd=\tempcounta                       
\vertadjust[\tokenc`\tokend`{#7}]
\ifnum\tempcounta<\tempcountd                 
\tempcounta=\tempcountd\fi                    
\advance \yext by\tempcounta                  
\advance \yoff by-\tempcounta                 %
\multiply \xext by2 
\leftsladjust[\tokenb`\tokena`{#3}]
\rightsladjust[\tokend`\tokena`{#5}]%
\begin{picture}(\xext,\yext)(\xoff,\yoff)%
\putAtrianglepairp
<\arrowtypea`\arrowtypeb`\arrowtypec`\arrowtyped`\arrowtypee;\height>%
(0,0)[#2;#3`#4`#5`#6`{#7}]%
\end{picture}%
}}

\def\putVtrianglep<#1>(#2,#3)[#4`#5`#6;#7`#8`#9]{{%
\settriparms[#1]%
\xpos=#2 \ypos=#3
\advance\ypos by\height
{\multiply\height by2
\puthmorphism(\xpos,\ypos)[#4`#5`{#7}]{\height}{\arrowtypea}a}%
\putmorphism(\xpos,\ypos)(1,-1)[`#6`{#8}]{\height}{\arrowtypeb}l%
\advance\xpos by\height
\advance\xpos by\height
\putmorphism(\xpos,\ypos)(-1,-1)[``{#9}]{\height}{\arrowtypec}r%
}}

\def\putVtriangle{\@ifnextchar <{\putVtrianglep}{\putVtrianglep
   <\arrowtypea`\arrowtypeb`\arrowtypec;\height>}}
\def\Vtriangle{\@ifnextchar <{\Vtrianglep}{\Vtrianglep
   <\arrowtypea`\arrowtypeb`\arrowtypec;\height>}}

\def\Vtrianglep<#1>[#2`#3`#4;#5`#6`#7]{{
\settriparms[#1]
\width=\height                         
\xext=\width                           
\yext=\height                          
\topadjust[#2`#3`{#5}]
\botadjust[#4``]
\multiply \xext by2 
\leftsladjust[#2`#3`{#6}]
\rightsladjust[#3`#4`{#7}]
\begin{picture}(\xext,\yext)(\xoff,\yoff)%
\putVtrianglep<\arrowtypea`\arrowtypeb`\arrowtypec;\height>%
(0,0)[#2`#3`#4;#5`#6`{#7}]%
\end{picture}%
}}

\def\putVtrianglepairp<#1>(#2)[#3;#4`#5`#6`#7`#8]{{
\settripairparms[#1]%
\setpos(#2)%
\settokens[#3]%
\advance\ypos by\height
\putmorphism(\xpos,\ypos)(1,-1)[`\tokend`{#6}]{\height}{\arrowtypec}l%
\puthmorphism(\xpos,\ypos)[\tokena`\tokenb`{#4}]{\height}{\arrowtypea}a%
\advance\xpos by\height
\putvmorphism(\xpos,\ypos)[``{#7}]{\height}{\arrowtyped}m%
\advance\xpos by\height
\putmorphism(\xpos,\ypos)(-1,-1)[``{#8}]{\height}{\arrowtypee}r%
}}

\def\putVtrianglepair{\@ifnextchar <{\putVtrianglepairp}{\putVtrianglepairp%
    <\arrowtypea`\arrowtypeb`\arrowtypec`\arrowtyped`\arrowtypee;\height>}}
\def\Vtrianglepair{\@ifnextchar <{\Vtrianglepairp}{\Vtrianglepairp%
    <\arrowtypea`\arrowtypeb`\arrowtypec`\arrowtyped`\arrowtypee;\height>}}

\def\Vtrianglepairp<#1>[#2;#3`#4`#5`#6`#7]{{%
\settripairparms[#1]%
\settokens[#2]
\xext=\height                  
\width=\height                 
\yext=\height                  
\vertadjust[\tokena`\tokenb`{#4}]
\tempcountd=\tempcounta        
\vertadjust[\tokenb`\tokenc`{#5}]
\ifnum\tempcounta<\tempcountd%
\tempcounta=\tempcountd\fi
\advance \yext by\tempcounta
\botadjust[\tokend``]%
\multiply \xext by2
\leftsladjust[\tokena`\tokend`{#6}]%
\rightsladjust[\tokenc`\tokend`{#7}]%
\begin{picture}(\xext,\yext)(\xoff,\yoff)%
\putVtrianglepairp
<\arrowtypea`\arrowtypeb`\arrowtypec`\arrowtyped`\arrowtypee;\height>%
(0,0)[#2;#3`#4`#5`#6`{#7}]%
\end{picture}%
}}

\def\putCtrianglep<#1>(#2,#3)[#4`#5`#6;#7`#8`#9]{{%
\settriparms[#1]%
\xpos=#2 \ypos=#3
\advance\ypos by\height
\putmorphism(\xpos,\ypos)(1,-1)[``{#9}]{\height}{\arrowtypec}l%
\advance\xpos by\height
\advance\ypos by\height
\putmorphism(\xpos,\ypos)(-1,-1)[#4`#5`{#7}]{\height}{\arrowtypea}l%
{\multiply\height by 2
\putvmorphism(\xpos,\ypos)[`#6`{#8}]{\height}{\arrowtypeb}r}%
}}

\def\putCtriangle{\@ifnextchar <{\putCtrianglep}{\putCtrianglep
    <\arrowtypea`\arrowtypeb`\arrowtypec;\height>}}
\def\Ctriangle{\@ifnextchar <{\Ctrianglep}{\Ctrianglep
    <\arrowtypea`\arrowtypeb`\arrowtypec;\height>}}

\def\Ctrianglep<#1>[#2`#3`#4;#5`#6`#7]{{
\settriparms[#1]
\width=\height                          
\xext=\width                            
\yext=\height                           
\multiply \yext by2 
\topadjust[#2``]
\botadjust[#4``]
\sladjust[#3`#2`{#5}]{\width}
\tempcountd=\tempcounta                 
\sladjust[#3`#4`{#7}]{\width}
\ifnum \tempcounta<\tempcountd          
\tempcounta=\tempcountd\fi              
\advance \xext by\tempcounta            
\advance \xoff by-\tempcounta           %
\rightadjust[#2`#4`{#6}]
\begin{picture}(\xext,\yext)(\xoff,\yoff)%
\putCtrianglep<\arrowtypea`\arrowtypeb`\arrowtypec;\height>%
(0,0)[#2`#3`#4;#5`#6`{#7}]%
\end{picture}%
}}

\def\putDtrianglep<#1>(#2,#3)[#4`#5`#6;#7`#8`#9]{{%
\settriparms[#1]%
\xpos=#2 \ypos=#3
\advance\xpos by\height \advance\ypos by\height
\putmorphism(\xpos,\ypos)(-1,-1)[``{#9}]{\height}{\arrowtypec}r%
\advance\xpos by-\height \advance\ypos by\height
\putmorphism(\xpos,\ypos)(1,-1)[`#5`{#8}]{\height}{\arrowtypeb}r%
{\multiply\height by 2
\putvmorphism(\xpos,\ypos)[#4`#6`{#7}]{\height}{\arrowtypea}l}%
}}

\def\putDtriangle{\@ifnextchar <{\putDtrianglep}{\putDtrianglep
    <\arrowtypea`\arrowtypeb`\arrowtypec;\height>}}
\def\Dtriangle{\@ifnextchar <{\Dtrianglep}{\Dtrianglep
   <\arrowtypea`\arrowtypeb`\arrowtypec;\height>}}

\def\Dtrianglep<#1>[#2`#3`#4;#5`#6`#7]{{
\settriparms[#1]
\width=\height                         
\xext=\height                          
\yext=\height                          
\multiply \yext by2 
\topadjust[#2``]
\botadjust[#4``]
\leftadjust[#2`#4`{#5}]
\sladjust[#3`#2`{#5}]{\height}
\tempcountd=\tempcountd                
\sladjust[#3`#4`{#7}]{\height}
\ifnum \tempcounta<\tempcountd         
\tempcounta=\tempcountd\fi             
\advance \xext by\tempcounta           %
\begin{picture}(\xext,\yext)(\xoff,\yoff)
\putDtrianglep<\arrowtypea`\arrowtypeb`\arrowtypec;\height>%
(0,0)[#2`#3`#4;#5`#6`{#7}]%
\end{picture}%
}}

\def\setrecparms[#1`#2]{\width=#1 \height=#2}%
\def\recursep<#1`#2>[#3;#4`#5`#6`#7`#8]{{%
\width=#1 \height=#2
\settokens[#3]
\settowidth{\tempdimen}{$\tokena$}
\ifdim\tempdimen=0pt
  \savebox{\tempboxa}{\hbox{$\tokenb$}}%
  \savebox{\tempboxb}{\hbox{$\tokend$}}%
  \savebox{\tempboxc}{\hbox{$#6$}}%
\else
  \savebox{\tempboxa}{\hbox{$\hbox{$\tokena$}\times\hbox{$\tokenb$}$}}%
  \savebox{\tempboxb}{\hbox{$\hbox{$\tokena$}\times\hbox{$\tokend$}$}}%
  \savebox{\tempboxc}{\hbox{$\hbox{$\tokena$}\times\hbox{$#6$}$}}%
\fi
\ypos=\height
\divide\ypos by 2
\xpos=\ypos
\advance\xpos by \width
\xext=\xpos \yext=\height
\topadjust[#3`\usebox{\tempboxa}`{#4}]%
\botadjust[#5`\usebox{\tempboxb}`{#8}]%
\sladjust[\tokenc`\tokenb`{#5}]{\ypos}%
\tempcountd=\tempcounta
\sladjust[\tokenc`\tokend`{#5}]{\ypos}%
\ifnum \tempcounta<\tempcountd
\tempcounta=\tempcountd\fi
\advance \xext by\tempcounta
\advance \xoff by-\tempcounta
\rightadjust[\usebox{\tempboxa}`\usebox{\tempboxb}`\usebox{\tempboxc}]%
\bfig
\putCtrianglep<-1`1`1;\ypos>(0,0)[`\tokenc`;#5`#6`{#7}]%
\puthmorphism(\ypos,0)[\tokend`\usebox{\tempboxb}`{#8}]{\width}{-1}b%
\puthmorphism(\ypos,\height)[\tokenb`\usebox{\tempboxa}`{#4}]{\width}{-1}a%
\advance\ypos by \width
\putvmorphism(\ypos,\height)[``\usebox{\tempboxc}]{\height}1r%
\efig
}}

\def\recurse{\@ifnextchar <{\recursep}{\recursep<\width`\height>}}

\def\puttwohmorphisms(#1,#2)[#3`#4;#5`#6]#7#8#9{{%
\puthmorphism(#1,#2)[#3`#4`]{#7}0a
\ypos=#2
\advance\ypos by 20
\puthmorphism(#1,\ypos)[\phantom{#3}`\phantom{#4}`#5]{#7}{#8}a
\advance\ypos by -40
\puthmorphism(#1,\ypos)[\phantom{#3}`\phantom{#4}`#6]{#7}{#9}b
}}

\def\puttwovmorphisms(#1,#2)[#3`#4;#5`#6]#7#8#9{{%
\putvmorphism(#1,#2)[#3`#4`]{#7}0a
\xpos=#1
\advance\xpos by -20
\putvmorphism(\xpos,#2)[\phantom{#3}`\phantom{#4}`#5]{#7}{#8}l
\advance\xpos by 40
\putvmorphism(\xpos,#2)[\phantom{#3}`\phantom{#4}`#6]{#7}{#9}r
}}

\def\puthcoequalizer(#1)[#2`#3`#4;#5`#6`#7]#8#9{{%
\setpos(#1)%
\puttwohmorphisms(\xpos,\ypos)[#2`#3;#5`#6]{#8}11%
\advance\xpos by #8
\puthmorphism(\xpos,\ypos)[\phantom{#3}`#4`#7]{#8}1{#9}
}}

\def\putvcoequalizer(#1)[#2`#3`#4;#5`#6`#7]#8#9{{%
\setpos(#1)%
\puttwovmorphisms(\xpos,\ypos)[#2`#3;#5`#6]{#8}11%
\advance\ypos by -#8
\putvmorphism(\xpos,\ypos)[\phantom{#3}`#4`#7]{#8}1{#9}
}}

\def\putthreehmorphisms(#1)[#2`#3;#4`#5`#6]#7(#8)#9{{%
\setpos(#1) \settypes(#8)
\if a#9 %
     \vertsize{\tempcounta}{#5}%
     \vertsize{\tempcountb}{#6}%
     \ifnum \tempcounta<\tempcountb \tempcounta=\tempcountb \fi
\else
     \vertsize{\tempcounta}{#4}%
     \vertsize{\tempcountb}{#5}%
     \ifnum \tempcounta<\tempcountb \tempcounta=\tempcountb \fi
\fi
\advance \tempcounta by 60
\puthmorphism(\xpos,\ypos)[#2`#3`#5]{#7}{\arrowtypeb}{#9}
\advance\ypos by \tempcounta
\puthmorphism(\xpos,\ypos)[\phantom{#2}`\phantom{#3}`#4]{#7}{\arrowtypea}{#9}
\advance\ypos by -\tempcounta \advance\ypos by -\tempcounta
\puthmorphism(\xpos,\ypos)[\phantom{#2}`\phantom{#3}`#6]{#7}{\arrowtypec}{#9}
}}

\def\putarc(#1,#2)[#3`#4`#5]#6#7#8{{%
\xpos #1
\ypos #2
\width #6
\arrowlength #6
\putbox(\xpos,\ypos){#3\vphantom{#4}}%
{\advance \xpos by\arrowlength
\putbox(\xpos,\ypos){\vphantom{#3}#4}}%
\horsize{\tempcounta}{#3}%
\horsize{\tempcountb}{#4}%
\divide \tempcounta by2
\divide \tempcountb by2
\advance \tempcounta by30
\advance \tempcountb by30
\advance \xpos by\tempcounta
\advance \arrowlength by-\tempcounta
\advance \arrowlength by-\tempcountb
\halflength=\arrowlength \divide\halflength by 2
\divide\arrowlength by 5
\put(\xpos,\ypos){\bezier{\arrowlength}(0,0)(50,50)(\halflength,50)}
\ifnum #7=-1 \put(\xpos,\ypos){\vector(-3,-2)0} \fi
\advance\xpos by \halflength
\put(\xpos,\ypos){\xpos=\halflength \advance\xpos by -50
   \bezier{\arrowlength}(0,50)(\xpos,50)(\halflength,0)}
\ifnum #7=1 {\advance \xpos by
   \halflength \put(\xpos,\ypos){\vector(3,-2)0}} \fi
\advance\ypos by 50
\vertsize{\tempcounta}{#5}%
\divide\tempcounta by2
\advance \tempcounta by20
\if a#8 %
   \advance \ypos by\tempcounta
   \putbox(\xpos,\ypos){#5}%
\else
   \advance \ypos by-\tempcounta
   \putbox(\xpos,\ypos){#5}%
\fi
}}

\makeatother

\usepackage{amsthm}
\usepackage{dsfont}
\usepackage{stmaryrd}

\newtheorem{theorem}{Theorem}[section]
\newtheorem{lemma}[theorem]{Lemma}

\newtheorem{proposition}[theorem]{Proposition}

\begin{document}

\sloppy

\newcommand{\nl}{\hspace{2cm}\\ }

\def\nec{\Box}
\def\pos{\Diamond}
\def\diam{{\tiny\Diamond}}

\def\lc{\lceil}
\def\rc{\rceil}
\def\lf{\lfloor}
\def\rf{\rfloor}
\def\lk{\langle}
\def\rk{\rangle}
\def\blk{\dot{\langle\!\!\langle}}
\def\brk{\dot{\rangle\!\!\rangle}}

\newcommand{\pa}{\parallel}
\newcommand{\lra}{\longrightarrow}
\newcommand{\hra}{\hookrightarrow}
\newcommand{\hla}{\hookleftarrow}
\newcommand{\ra}{\rightarrow}
\newcommand{\la}{\leftarrow}
\newcommand{\lla}{\longleftarrow}
\newcommand{\da}{\downarrow}
\newcommand{\ua}{\uparrow}
\newcommand{\dA}{\downarrow\!\!\!^\bullet}
\newcommand{\uA}{\uparrow\!\!\!_\bullet}
\newcommand{\Da}{\Downarrow}
\newcommand{\DA}{\Downarrow\!\!\!^\bullet}
\newcommand{\UA}{\Uparrow\!\!\!_\bullet}
\newcommand{\Ua}{\Uparrow}
\newcommand{\Lra}{\Longrightarrow}
\newcommand{\Ra}{\Rightarrow}
\newcommand{\Lla}{\Longleftarrow}
\newcommand{\La}{\Leftarrow}
\newcommand{\nperp}{\perp\!\!\!\!\!\setminus\;\;}
\newcommand{\pq}{\preceq}

\newcommand{\lms}{\longmapsto}
\newcommand{\ms}{\mapsto}
\newcommand{\subseteqnot}{\subseteq\hskip-4 mm_\not\hskip3 mm}

\def\o{{\omega}}

\def\bA{{\bf A}}
\def\bEM{{\bf EM}}
\def\bM{{\bf M}}
\def\bN{{\bf N}}
\def\bF{{\bf F}}
\def\bC{{\bf C}}
\def\bI{{\bf I}}
\def\bK{{\bf K}}
\def\bL{{\bf L}}
\def\bT{{\bf T}}
\def\bS{{\bf S}}
\def\bD{{\bf D}}
\def\bB{{\bf B}}
\def\bW{{\bf W}}
\def\bP{{\bf P}}
\def\bX{{\bf X}}
\def\bY{{\bf Y}}
\def\ba{{\bf a}}
\def\bb{{\bf b}}
\def\bc{{\bf c}}
\def\bd{{\bf d}}
\def\bh{{\bf h}}
\def\bi{{\bf i}}
\def\bj{{\bf j}}
\def\bk{{\bf k}}
\def\bm{{\bf m}}
\def\bn{{\bf n}}
\def\bp{{\bf p}}
\def\bq{{\bf q}}
\def\be{{\bf e}}
\def\br{{\bf r}}
\def\bi{{\bf i}}
\def\bs{{\bf s}}
\def\bt{{\bf t}}
\def\jeden{{\bf 1}}
\def\dwa{{\bf 2}}
\def\trzy{{\bf 3}}

\def\cB{{\cal B}}
\def\cA{{\cal A}}
\def\cC{{\cal C}}
\def\cD{{\cal D}}
\def\cE{{\cal E}}
\def\cEM{{\cal EM}}
\def\cF{{\cal F}}
\def\cG{{\cal G}}
\def\cI{{\cal I}}
\def\cJ{{\cal J}}
\def\cK{{\cal K}}
\def\cL{{\cal L}}
\def\cN{{\cal N}}
\def\cM{{\cal M}}
\def\cO{{\cal O}}
\def\cP{{\cal P}}
\def\cQ{{\cal Q}}
\def\cR{{\cal R}}
\def\cS{{\cal S}}
\def\cT{{\cal T}}
\def\cU{{\cal U}}
\def\cV{{\cal V}}
\def\cW{{\cal W}}
\def\cX{{\cal X}}
\def\cY{{\cal Y}}


\def\Mnd{{\bf Mnd}}
\def\AMnd{{\bf AnMnd}}
\def\An{{\bf An}}
\def\Poly{{\bf Poly}}
\def\San{{\bf San}}
\def\cSan{{\bf cSan}}
\def\Taut{{\bf Taut}}
\def\PMnd{{\bf PolyMnd}}
\def\SanMnd{{\bf SanMnd}}
\def\RiMnd{{\bf RiMnd}}
\def\End{{\bf End}}
\def\cEnd{{\bf cEnd}}

\def\ET{\bf ET}
\def\RegET{\bf RegET}
\def\RET{\bf RegET}
\def\LrET{\bf LrET}
\def\RiET{\bf RiET}
\def\SregET{\bf SregET}
\def\Cart{\bf Cart}
\def\wCart{\bf wCart}
\def\CartMnd{\bf CartMnd}
\def\wCartMnd{\bf wCartMnd}

\def\LT{\bf LT}
\def\RegLT{\bf RegLT}
\def\ALT{\bf AnLT}
\def\RiLT{\bf RiLT}

\def\FOp{\bf FOp}
\def\RegOp{\bf RegOp}
\def\SOp{\bf SOp}
\def\RiOp{\bf RiOp}

\def\bCat{{{\bf Cat}}}
\def\MonCat{{{\bf MonCat}}}
\def\Mon{{{\bf Mon}}}
\def\Cat{{{\bf Cat}}}

\def\F{\mathds{F}}
\def\S{\mathds{S}}
\def\I{\mathds{I}}
\def\B{\mathds{B}}

\def\V{\mathds{V}}
\def\W{\mathds{W}}
\def\M{\mathds{M}}
\def\N{\mathds{N}}
\def\R{\mathds{R}}

\def\Op{{\cal O}p}

\def\Vb{\bar{\mathds{V}}}
\def\Wb{\bar{\mathds{W}}}

\def\Sym{{\cal S}ym}

\pagenumbering{arabic} \setcounter{page}{1}

\title{\bf\Large co-Semi-analytic functors}

\author{ Marek Zawadowski\\
}

\maketitle
\begin{abstract} We characterize the category of co-semi-analytic functors and describe an action of semi-analytic functors on co-semi-analytic functors.
\end{abstract}


\section{Introduction} Let $\I$ be the (skeleton of the) category of finite sets and monomorphisms.
The functor $Set^{\I^{op}}\ra Cat(Set^op,Set)$ of the left K extension along $\I^{op}\ra Set^{op}$ is conservative. The essential image of this functor is a category of functors that are in a sense a dual presentation to the semi-analytic functors \cite{SZ}. This is why we call this category the category of  co-semi-analytic functors. In this note, we give an abstract characterization of this category. Moreover, we show that the category of semi-analytic functors acts on this category and we show some examples of the actions along this action.

Contravariant functors on finite sets were considered in \cite{P}, \cite{D}.

\subsection*{Notation}
Let $[n]=\{0,\ldots, n\}$, $(n]=\{1,\ldots, n\}$, $\o$ - denote the set of natural numbers. The set $X^n$, $n$-th power of $X$, is interpreted as $X^{(n]}$, a set of unctions, when convenient. The skeletal category equivalent to the category of finite sets $Set_{fin}$ will be denoted by $\mathds{F}$. We will be assuming that the objects of $\mathds{F}$ are sets $(n]$, for $n\in \o$.  The subcategories of $\mathds{F}$ with the same objects as $\mathds{F}$ but having as objects bijection, surjections and injections will be denoted by $\mathds{B}$, $\mathds{S}$, $\mathds{I}$, respectively. $S_n$ is the group of permutations of $(n]$. When $S_n$ acts on a set $A$ on the right and on the set $B$ on the left, then the set $A\otimes_nB$ is the usual tensor product of $S_n$-sets. Let $Epi(X,(n])$ denote the set of epimorphisms from the set $X$ to $(n]$. $S_n$ acts on $Epi(X,(n])$ on the left by compositions. If $A:\I^{op}\ra Set$ is a functor, $f:(n]\ra (m]$, $a\in A_m\;(=A(m])$,  then we often write $a\cdot_A f$ instead of $A(f)(a)$. $S_n$ acts of the right on $A_n$ and, according with the previous notation, we write $a\cdot_A \sigma$ instead of $A(\sigma)(a)$. Thus we can form a set
\[ A_n\otimes_n Epi(X,(n]) \]
whose elements are equivalence classes of pairs $\lk a,\overleftarrow{x}\rk$ such that $a\in A_n$ and $\overleftarrow{x}:X\ra (n]$ epi\footnote{To emphasize the domain of an epi $X\ra(n]$, we name it $\overleftarrow{x}$ as it is in a sense a dual to $(n]\ra X$. And it is natural to name the later map $\overrightarrow{x}$.}. We identify pairs
\[ \lk a,\sigma\circ \overleftarrow{x}\rk \sim  \lk a\cdot\sigma, \overleftarrow{x}\rk \]
for $\sigma\in S_n$.

\section{The category of co-semi-analytic functors}\label{co-semianalytic}

A natural transformation $\tau : F\ra G:\cC\ra \cD$ is {\em semi-cartesian} iff the naturality squares for monomorphisms in $\cC$ are pullbacks in $\cD$. Let $\cEnd$ denote the category $Nat(Set^{op},Set)$, i.e. the category of contravariant functors on $Set$ and natural transformations.

We define a functor
\[ \check{(-)} : Set^{\I^{op}} \lra \cEnd \]
Let $A: \I^{op}\ra Set$. We put
\[ \check{A}(X) = \sum_{n\in \o} A_n\otimes_n Epi(X,(n]) \]
For a function $f:Y\ra X$ and $[a,\overleftarrow{x}]_\sim \in  A_n\otimes_n Epi(X,(n])$ we put
\[ \check{A}(f)([a,\overleftarrow{x}]_\sim) = [a\cdot_A{f'},\overleftarrow{x}'] \]
where the square
\begin{center} \xext=600 \yext=500
\begin{picture}(\xext,\yext)(\xoff,\yoff)
\setsqparms[3`-1`-2`3;500`400]
 \putsquare(0,50)[X`(n]`Y`{(m]};\overleftarrow{x}`f`f'`\overleftarrow{x}']
 \end{picture}
\end{center}
commutes and $\overleftarrow{x}'$, $f'$ is the epi-mono factorization of $\overleftarrow{x}\circ f$.

Let $\tau:A\ra B$ be a natural transformation in $Set^{\I^{op}}$. For $[a,\overleftarrow{x}]_\sim \in  A_n\otimes_n Epi(X,(n]) \subseteq \check{A}(X)$ we put
\[ \check{\tau}_X([a,\overleftarrow{x}]_\sim) = [\tau_n(a),\overleftarrow{x}]_\sim \]

We have
\begin{proposition}\label{check vs kan}
The functor $\check{(-)}: Set^{\I^{op}} \lra \cEnd$ is well defined, and it is isomorphic to the left Kan extension along the inclusion functor $\I^{op}\ra Set^{op}$.
\end{proposition}

{\em Proof.} The fact that $\check{(-)}$ is well defined is easy. It is well known that the left Kan extension can be calculated with coends. It is also easy to check that, for $A\in Set^{\I^{op}}$ and a set $X$, we have the second isomorphism
\[ \check{A}(X) = \int^{(n]\in \I^{op}}A_n\times Set(X,(n]) \cong \sum_{n\in\o} A_n\otimes_n Epi(X,(n]) \]
\[ \hskip 3cm [a,f]_\sim \mapsto [a\cdot f',q] \]
where $f=f'\circ q$ is the epi-mono factorization of $f$. The details are left for the reader.
$\boxempty$
\vskip 2mm

The {\em canonical finitary cone $\gamma$ under a set $X$} in $Set$ is the cone from the vertex $X$ to the functor
\[ X\da \F \stackrel{\pi_X}{\lra} \F \lra Set \]
\[ X\ra (n] \mapsto (n] \]
such that $\gamma_f=f$, for the $f:X\ra (n]$ in $X\da\F$. The {\em canonical finitary cocone $\kappa$ over a set $X$} in $Set^{op}$ is the dual of the canonical finitary cone under $X$ i.e. the cocone from the functor
\[ (X\da \F)^{op} \stackrel{(\pi_X)^{op}}{\lra} \F^{op} \lra Set^{op} \]
to the vertex $X$, such that $\kappa_f=f$ for $f:X\ra (n]$ in $(X\da\F)^{op}$.

\begin{theorem}\label{T-char-csan}
The functors $\check{(-)}: Set^{\I^{op}} \lra \cEnd$ is conservative and its essential image consists of functors sending
\begin{enumerate}
  \item pullbacks along monos to pullbacks;
  \item canonical finitary cocone under $X$ in $Set^{op}$ to a colimiting cocone in $Set$, for any set $X$;
\end{enumerate}
and semi-analytic natural transformations.
\end{theorem}

Because of the above theorem, by definition, the essential image of the functor $\check{(-)}$ is called the category of co-semi-analytic functors, and is denoted $\cSan$.

We shall prove the above theorem through a series of Lemmas.

\begin{lemma} \label{L-im}Let $\tau:A\ra B$ be a morphism in $Set^{\I^{op}}$. Then
\begin{enumerate}
  \item the functor $\check{A} :Set^{op}\ra Set$ sends
  \begin{enumerate}
  \item pullbacks along monos to pulbacks;
  \item canonical finitary cocone over $X$ in $Set^{op}$ to a colimiting cocone in $Set$, for any set $X$;
  \end{enumerate}
  \item the natural transformation $\check{\tau}: \check{A}\ra \check{B}$ is semi-analytic.
\end{enumerate}
\end{lemma}

{\em Proof.} Let the square
\begin{center} \xext=800 \yext=650
\begin{picture}(\xext,\yext)(\xoff,\yoff)
\setsqparms[1`1`1`1;800`450]
 \putsquare(0,100)[X`Z`Y`T;g`f`f'`g']
 \end{picture}
\end{center}
be a pushout in $Set$ with $f$ (and hence also $f'$) epi. We need to show that the square
\begin{center} \xext=800 \yext=650
\begin{picture}(\xext,\yext)(\xoff,\yoff)
\setsqparms[-1`-1`-1`-1;800`450]
 \putsquare(0,100)[\check{A}(X)`\check{A}(Z)`\check{A}(Y)`\check{A}(T);\check{A}(g)`\check{A}(f)`\check{A}(f')`\check{A}(g')]
 \end{picture}
\end{center}
is a pullback in $Set$. Let $[a,\overleftarrow{y}: Y\ra (n]]_\sim \in \check{A}(Y)$ and $[b,\overleftarrow{z}: Z\ra (m]]_\sim \in \check{A}(Z)$ such that
\[ \check{A}(f)([a,\overleftarrow{y}: Y\ra (n]]_\sim) = \check{A}(g)([b,\overleftarrow{z}: Z\ra (m]]_\sim) \]
Thus we have a mono $\bar{f}:(n]\ra (m]$ such that
\[ \bar{f}\circ \overleftarrow{y}\circ f = \overleftarrow{z}\circ g \;\;\; {\rm and } \;\;\; b\cdot_A\bar{f}= a  \]
Hence there is  a unique morphisms $\overleftarrow{t} : T \ra (m]$ such that
\[ \overleftarrow{t}\circ g'= \bar{f}\circ \overleftarrow{y}\;\;\;{\rm and }\;\;\; \overleftarrow{t}\circ f' = \overleftarrow{z} \]
As $\overleftarrow{z}$ is epi, $\overleftarrow{t}$ is epi, as well. Then $[b, \overleftarrow{t}]_\sim \in \check{A}(T)$ and we have
\[ \check{A}(f')([b,\overleftarrow{t}]_\sim )= [b,\overleftarrow{t}\circ f']_\sim  =  [b,\overleftarrow{z}]_\sim   \]

\[ \check{A}(g')([b,\overleftarrow{t}]_\sim ) = [b\cdot_A\bar{f},\overleftarrow{y}]_\sim = [a,\overleftarrow{y}]_\sim  \]
As any functor $Set^{op}\ra Set$ sends monos to monos, it follows that $\check{A}$ sends pullbacks along monos to pullbacks.

Let $X$ be a set. We shall show that $\check{A}$ sends the canonical finitary cocone over $X$ in $Set^{op}$ to a colimiting cocone in $Set$. Let $[a,\overleftarrow{x}: X\ra (n]]_\sim \in \check{A}(X)$.
Then $[a,1_n]_\sim \in \check{A}(n]$ and
\[ \check{A}(\overleftarrow{x})([a,1_n]_\sim) = [a,\overleftarrow{x}]_\sim \]
Thus the morphisms $\check{A}(\overleftarrow{x}):\check{A}(n]\ra \check{A}(X)$ with $\overleftarrow{x}:X\ra (n]$ epi, jointly cover $\check{A}(X)$.

If we have another element $[b,q:(m]\ra (k]]_\sim\in \check{A}(m]$ and $f:X\ra (m]$ a function such that
\[ \check{A}(f)([b,q]_\sim)=[a,\overleftarrow{x}]_\sim\]
then with $q'$, $f'$ being the epi-mono factorization (see diagram below) of $q\circ f$ we have
\[  \check{A}(f)([b,q]_\sim)=[b\cdot_A f',q']_\sim = [a,\overleftarrow{x}]_\sim \]
Hence we have a $\sigma\in S_n$ making the left triangle
\begin{center} \xext=800 \yext=700
\begin{picture}(\xext,\yext)(\xoff,\yoff)
\setsqparms[1`1`1`1;500`500]
 \putsquare(500,50)[X`(m]`{(n]}`{(k]};f`q'`q`f']
 \putmorphism(0,50)(1,0)[(n]`\phantom{(n]}`\sigma]{500}{1}a
 \put(100,300){\makebox(100,100){$\overleftarrow{x}$}}
 \put(410,530){\vector(-1,-1){400}}
 \end{picture}
\end{center}
commute and \[ a=(b\cdot_A f')\cdot_A \sigma\]
Thus we have a commuting square
\begin{center} \xext=800 \yext=700
\begin{picture}(\xext,\yext)(\xoff,\yoff)
\setsqparms[-1`-1`-1`-1;800`500]
 \putsquare(0,50)[\check{A}(X)`\check{A}(m]`\check{A}(n]`{\check{A}(k]};\check{A}(f)`\check{A}(\overleftarrow{x})`\check{A}(q)`\check{A}(f'\circ\sigma)]
 \end{picture}
\end{center}
and $[b,1_k]_\sim\in \check{A}(k]$ such that
\[ \check{A}(q) ([b,1_k]_\sim) =[b,q]_\sim \]
and
\[ \check{A}(f'\circ\sigma) ([b,1_k]_\sim) =[a,1_n]_\sim \]
i.e. if two elements go to the same element in $\check{A}(X)$, they are related in the cocone. Thus $\check{A} :Set^{op}\ra Set$ sends canonical finitary cocone under $X$ in $Set^{op}$ to a colimiting cocone in $Set$, for any set $X$.

It remains to show that the natural transformation $\check{\tau}:\check{A}\ra\check{B}$ is semi-cartesian. Let $g:X\ra Y$ be an epi in $Set$. We shall show that the square
\begin{center} \xext=800 \yext=650
\begin{picture}(\xext,\yext)(\xoff,\yoff)
 \setsqparms[1`-1`-1`1;800`450]
 \putsquare(0,100)[\check{A}(X)`\check{B}(X)`\check{A}(Y)`\check{B}(Y);\check{\tau}_X`\check{A}(g)`\check{B}(g)`\check{\tau}_Y]
\end{picture}
\end{center}
is a pullback. Fix $[a,\overleftarrow{x}: X\ra (n]]_\sim \in \check{A}(X)$ and $[b,\overleftarrow{y}: Y\ra (n]]_\sim \in \check{B}(Y)$ such that
\[ [\tau_n(a),\overleftarrow{x}]_\sim= \check{\tau}([a,\overleftarrow{x}]_\sim) =  \check{B}(g)([b,\overleftarrow{y}]_\sim) =[b,\overleftarrow{y}\circ g]_\sim  \]
Thus we have a permutation $\sigma : (n]\ra (n]$ such that
\[ \tau_n(a)\cdot_B\sigma= b  \;\;\; {\rm and } \;\;\; \overleftarrow{x}=\sigma\circ \overleftarrow{y}\circ g  \]
Then $[a,\sigma\circ \overleftarrow{y}]_\sim \in \check{A}(Y)$ and we have
\[ \check{A}(g)([a,\sigma\circ \overleftarrow{y}]_\sim) = [a,\sigma\circ \overleftarrow{y}\circ g]_\sim = [a,\overleftarrow{x}]_\sim \]

\[ \check{\tau}_Y([a,\sigma\circ \overleftarrow{y}]_\sim) = [\tau_n(a),\sigma\circ \overleftarrow{y}]_\sim = [\tau_n(a)\cdot_B\sigma, \overleftarrow{y}]_\sim = [b, \overleftarrow{y}]_\sim \]
As $\check{A}(g)$ is mono, the above square is a pullback. Thus $\check{\tau}$ is a semi-cartesian natural transformation. $\boxempty$
\vskip 2mm

\begin{lemma} \label{L-fu}Let $\cA :Set^{op}\ra Set$ be a functor that sends pullbacks along monos to pullbacks and canonical finitary cocone over $X$ in $Set^{op}$ to a colimiting cocone in $Set$, for any set $X$. Then there is a functor $A:\I^{op}\ra Set$ such that
$\check{A}$ is isomorphic to $\cA$.
\end{lemma}

{\em Proof.} Let $\cA :Set^{op}\ra Set$ be a functor with the properties described in Lemma. We define the functor $A:\I^{op}\ra Set$ as follows
\[ A_n= \cA(n] - \bigcup_{h} im (\cA(h)) \]
where the sum is over $0\leq m<n$ and (proper) epis $h:(n]\ra (m]$. The set $im (\cA(h)$ is the image $\cA(h)$ in $\cA(n]$. Then, for a mono $f:(n]\ra (m]$ in $\I$ we put
\[ A(f) =  \cA(f)_{\lceil A_m} \]
i.e. $A(f)$ is a restriction of $\cA(f)$ to a function $A_m\ra A_n$.

First we need to verify that $A$ is a well defined functor, i.e. that the restriction is the function with the appropriate domain and codomain.
Let $f:(n']\ra (n]$ be a mono $x\in \cA(n]$ but $\cA(f)(x)\not\in A_{n'}$. Thus there is a proper epi $h':(n]\ra (m']$ and $y\in \cA(m']$ such that
$\cA(h')(y)=\cA(f)(x)$. The square
\begin{center} \xext=800 \yext=650
\begin{picture}(\xext,\yext)(\xoff,\yoff)
 \setsqparms[-1`1`1`-1;800`450]
 \putsquare(0,100)[\cA(n]`\cA(m]`\cA(n']`{\cA(m']};\cA(h)`\cA(f)`\cA(f')`\cA(h')]
\end{picture}
\end{center}
is a pullback, where the square below
\begin{center} \xext=800 \yext=650
\begin{picture}(\xext,\yext)(\xoff,\yoff)
 \setsqparms[1`-1`-1`1;800`450]
 \putsquare(0,100)[(n]`(m]`{(n']}`{(m']};h`f`f'`h']
\end{picture}
\end{center}
is a pullback in $Set^{op}$, i.e. pushout of a proper epi $h'$ along $f$. Hence $h$ is also proper epi. Thus, there is $z\in \cA(m]$ such that
\[ \cA(h)(z)=x\;\;\; {\rm and}\;\;\; \cA(f')(z)=y \]
But this means that $x\not\in A_n$. This shows that $A(f)$ is well defined.

Now we define a natural isomorphism
\[ \varphi^\cA: \check{A} \ra \cA \]
so that
\[  \varphi^\cA_X([a,\overleftarrow{x}]_\sim ) = \cA(\overleftarrow{x})(a) \]
for any $X$, $n\in\o$, $a\in A_n$, and epi $\overleftarrow{x}:X\ra (n]$.

The fact that $\phi^\cA$ is a natural transformation is left for the reader.  $\phi^\cA$ is onto as $\cA$ sends canonical finitary cocones under any set $X$ in $Set^{op}$ to colimiting cocones in $Set$.

We shall show that $\phi^\cA$ is mono. Fix a set $X$ and let $[a,\overleftarrow{x}:X\ra (n]]_\sim, [a',\overleftarrow{x'}:X\ra (n']]_\sim\in \check{A}(X)$ such that
\[ \varphi^\cA_X([a,\overleftarrow{x}]_\sim )= \varphi^\cA_X([a',\overleftarrow{x'}]_\sim )  \]
By assumption, the pushout of epi $\overleftarrow{x}$ along epi $\overleftarrow{x'}$ in $Set$
\begin{center} \xext=800 \yext=650
\begin{picture}(\xext,\yext)(\xoff,\yoff)
 \setsqparms[-1`1`1`-1;800`450]
 \putsquare(0,100)[(n]`X`{(m]}`{(n']};\overleftarrow{x}`f`\overleftarrow{x'}`f']
\end{picture}
\end{center}
is sent to the pullback in $Set$ by $\cA$. Hence there is $a''\in \cA(m]$ such that
\[  \cA(f)(a'')=a  \;\;\; {\rm and}\;\;\; \cA(f')(a'')=a' \]
Thus, by definition of $A_n$ and $A_{n'}$ both $f$ and $f'$ are bijections (as they are epi but cannot be proper epi). Hence
\[ \cA(f^{-1}\circ f')(a)=a' \]
that is $[a,\overleftarrow{x}:X\ra (n]]_\sim = [a',\overleftarrow{x'}:X\ra (n]]_\sim$ and $\phi^\cA_X$ is mono, for any set $X$.
$\boxempty$
\vskip 2mm

\begin{lemma} \label{L-nt}Let $A,B :\I^{op}\ra Set$ be functors and $\psi:\check{A}\ra\check{B}$ a semi-analytic natural transformation. Then there is a natural transformation $\tau :A\ra B$ in $Set^{\I^{op}}$ such that $\check{\tau}=\psi$.
\end{lemma}
{\em Proof.} Let $\psi:\check{A}\ra\check{B}$ be a semi-analytic natural transformation. We define a natural transformation $\tau:A\ra B \in Set^{\I^{op}}$, as follows. Fix  $m\in\o$ and $a\in A_m$. Let $\psi_{(m]}([a,1_{(m]}]_\sim)=[b,p:(m]\ra(k]]_\sim\in \check{B}(m]$. As $p$ is an epi and $\psi$ is semi-analytic, the naturality square for $p$
\begin{center} \xext=800 \yext=650
\begin{picture}(\xext,\yext)(\xoff,\yoff)
 \setsqparms[1`-1`-1`1;800`450]
 \putsquare(0,100)[\check{A}(m]`\check{B}(m]`\check{A}(k]`{\check{B}(k]};\psi_{(m]}`\check{A}(p)`\check{B}(p)`\psi_{(k]}]
\end{picture}
\end{center}
is a pullback, and $\check{B}(p)([b,1_k]_\sim)=[b,p]_\sim$. Hence there is $[c,q]_\sim\in \check{A}(k]$ such that
\[  \check{A}(p)([c,q]_\sim)=[a,1_{(m]}]_\sim  \;\;\; {\rm and } \;\;\; \psi_{(k]}([c,q:(k]\ra(l]]_\sim)=[b,1_{(k]}]_\sim \]
In particular $q\circ p$ is a bijection, as $1_{(m]}$ is. Hence $k=m$ and $p$ is a bijection. We put
\[ \tau_m(a)= b\cdot_B p \]
Thus
\[ \psi_{(m]}([a,1_{(m]}]_\sim)=[\tau_m(a),1_{(m]}]_\sim \]
for $m\in\o$ and $a\in A_m$.
From the naturality of $\psi$ on a mono $f:(m']\ra (m]\in \I$, for any $a\in A_m$, we have
\[     [\tau_{m'}((a\cdot_A f),1_{(m']}]_\sim =\psi_{(m']}([a\cdot_Af),1_{(m']}]_\sim)= \]
\[ = \psi_{(m']}\circ \check{A}(f)([a,1_{(m]}]_\sim)  =\check{B}(f)\circ \psi_{(m]}([a,1_{(m]}]_\sim) =\]
\[= \check{B}(f)([\tau_m(a),1_{(m]}]_\sim) = [\tau_m(a)\cdot_B f,1_{(m]}]_\sim \]
Thus
\[ \tau_{m'}(A(f)(a))= B(f)(\tau_m(a)) \]
and hence $\tau: A\ra B$ is natural.

It remains to show that $\check{\tau}=\psi$. Fix set $X$, $m\in\o$ and $[a,q:X\ra(m]]_\sim \in \check{A}(X)$. Using naturality of $\psi$, $\check{\tau}$ on $q$ and the above, we have
\[ \psi_{X}([a,q]_\sim)=  \psi_{X}\circ \check{A}(q)([a,1_{(m]}]_\sim)=\]

\[=  \check{B}(q)\circ\psi_{(m]}([a,1_{(m]}]_\sim)   =\check{B}(q)([\tau_m(a),1_{(m]}]_\sim) =\]

\[  = \check{B}(q)\circ \check{\tau}_{(m]}([a,1_{(m]}]_\sim) = \check{\tau}_{X}\circ\check{A}(q) ([a,1_{(m]}]_\sim)  =\]

\[  =  \check{\tau}_{X}([a,q]_\sim) \]
Thus $\check{\tau}=\psi$, as required. $\boxempty$
\vskip 2mm

The fact that $\check{(-)}$ is immediate from the definition and hence Theorem \ref{T-char-csan} follows from Lemmas \ref{L-im}, \ref{L-fu}, \ref{L-nt}.

\section{The action}

The category $\End$ of endofunctors on $Set$ with composition as a tensor is strict monoidal. The composition
\[ \End \times  \cEnd \lra \cEnd \]
is an action of a monoidal category on a category. By Theorem \ref{T-char-csan} and characterization of semi-analytic functors (Theorem 2.2 of \cite{SZ}) the composition of semi-analytic functor with co-semi-analytic functor is co-semi-analytic. Thus the above action restricts to the action
\[ \San \times \cSan \lra \cSan \]
where $\San$ is the category of semi-analytic functors defined in \cite{SZ}. This category is equivalent to $Set^\S$ and has nice abstract characterization, see Section 2 of \cite{SZ}.

\vskip 2mm
{\bf Examples.}
The examples of actions of (semi-analytic) monads on $Set$ on contravariant functors on $Set$ consist of functors that build algebras of a (semi-analytic) monad out of sets in a contravariant way. We list some such examples below.
\begin{enumerate}
    \item Let $R$ be any algebra for a monad $T$ on $Set$ and let $\cR(X)$ be the algebra of functions on $X$ with values in $R$ with operations of the monad defined pointwise. Then $\cR$ is a contravariant functor on $Set$ on which we have an action of the monad $T$. The following two examples are, in a sense, special cases of this situation.
  \item The contravariant power-set functor $\cP : Set^{op}\ra Set$ is co-semi-analytic (see next example). Its value  on a set $X$ is the universe of a Boolean algebras $\cP(X)$ and the inverse image functor preserves the Boolean algebra structure. Thus the monad $T_{ba}$ on $Set$ for Boolean algebras acts on $\cP$. However, the monad $T_{ba}$ is not semi-analytic. 
  \item \label{ex} For any $n\in\o$ the functor
  \[ \cE^n : Set^{op}\lra Set \]
  \[ X\mapsto (n]^X \]
  is co-semi-analytic. The coefficient functor
  \[ E^n : \I^{op} \ra Set \]
  for $\cE^n$ is representable by $(n]$. The action of the monad $T_{ba}$ on the functor $\cP$ described above is a special case of the following. Let $T$ be a monad on $Set$ and $((n],\alpha:T(n]\ra (n])$ be a $T$-algebra. Then on $\cE^n(X)$ there is a natural structure of a $T$ algebra defined pointwise or using strength of $T$.

  Recall that if $X$ and $Y$ are sets and $x\in X$, then we have a function $\bar{x}: Y\ra X\times Y$ such that $\bar{x}(y)=\lk x,y\rk$. This allows to define strength on $T$:
   \[ st_{X,Y} : X\times T(Y)\lra T(X\times Y) \]
   so that  $st_{X,Y}(x,t)=\bar{x}(t)$.  It is an easy exercise to show that the strength on any semi-analytic monad is semi-cartesian.

   The $T$-algebra structure on $\cE^n(X)$ is the exponential adjoint of
  \begin{center} \xext=2800 \yext=150
\begin{picture}(\xext,\yext)(\xoff,\yoff)
  \putmorphism(0,0)(1,0)[X\times T((n]^X)`T(X\times (n]^X)`st_{X,(n]^X}]{1200}{1}a
  \putmorphism(1200,0)(1,0)[\phantom{T(X\times (n]^X)}`T(n]`T(ev)]{1000}{1}a
  \putmorphism(2200,0)(1,0)[\phantom{T(n]}`(n]`\alpha]{600}{1}a
\end{picture}
\end{center}
where $ev$ is the usual evaluation map.
\end{enumerate}



\end{document}